\input amstex
\documentstyle{amsppt}
\NoBlackBoxes
\magnification 1200
\vcorrection{-8mm}
\input epsf

\def\Q  {\Bbb Q}
\def\R  {\Bbb R}
\def\C  {\Bbb C}
\def\P  {\Bbb P}
\def\CP {\Bbb{CP}}
\def\baromega{\overline\omega}

\def\aa {\bold a}

\def\aaa{\widetilde\aa}

\let\wh=\widehat

\def\sectSym      {1}
\def\sectBirat      {1.1}
\def\sectConstr     {1.2}
\def\sectQP         {1.3}
\def\sectB          {1.4}
\def\sectPerturb  {2}
\def\sectDualHesse  {2.1}
\def\sectGenPerturb {2.2}
\def\sectStepOne    {2.3}
\def\sectStepTwo    {2.4}
\def\sectSeries     {2.5}

\def\propDegTen  {1}
\def\propPerturb {2}

\def\remGenLine {1}
\def\remMethod  {2}

\def\remUseSym  {3}
\def\remGLbis   {4}
\def\remHDbis   {5}

\def\refDHS  {1}
\def\refGLSa {2}
\def\refGLSb {3}

\def\refMO   {4}
\def\refS    {5}

\def\eqTenLines {1}
\def\eqLzero  {2}
\def\eqDet    {3}
\def\eqA      {4}
\def\eqPhi    {5}
\def\eqE      {6}
\def\eqSys    {7}

\def\figSym     {1}
\def\figCr      {2}
\def\figPerturb {3}

\topmatter
\title
       On curves of degree 10 with 12 triple points
\endtitle
\author
       S.~Yu.~Orevkov
\endauthor

\abstract
We construct an irreducible rational curve of degree 10 in $\Bbb{CP}^2$ which has 12 triple points, and
a union of three rational quartics with 19 triple points. This gives
counter-examples to a conjecture by Dimca, Harbourne, and Sticlaru.
We also prove that there exists an analytic family $C_u$ of curves of degree 10 with 12 triple points
which tends as $u\to 0$
to the union of the dual Hesse arrangement of lines (9 lines with 12 triple points) with an additional line.
We hope that our approach to the proof of the latter fact could be of independent interest.
\endabstract

\address IMT, Univ.~de Toulouse, Toulouse, France; Steklov Math.~Inst., Moscow, Russia.
\endaddress

\endtopmatter

\document
Recently Joaquim Ro\'e (private communication) asked if there exists
an irreducible curve of degree 10 in $\CP^2$ which has 12 ordinary triple points.
If such a curve exists, it should be rational by the genus formula.
We give an affirmative answer to that question.
This is a counterexample to [\refDHS, Conjecture 1.6], which states that {\sl there is no reduced plane curve
of degree $>9$ whose components are all rational and whose singularities (including infinitely near points)
all have multiplicity 3}
(note that this conjecture is presented in [\refDHS] as ``a more accessible target for a counterexample'').
Besides an irreducible curve of degree 10 with 12 triple points, we also construct a union of three rational
quartics with 19 triple points. Thus the bound 9 in the conjecture should be raised at least to 12.

\proclaim{ Proposition \propDegTen } (a).
The degree $10$ curve parameterized by $t\mapsto(x(t):y(t):z(t))$,
$$
 x = (t^3 + 2)(t^6 + 3 t^3 + 3), \quad  y = t(t^3+1)(t^3+2)(t^3+3), \quad z = t^9 + 3 t^6 - 3,
$$
has $12$ ordinary triple points.
Its Cartesian equation is $F(x,y,z)=0$ where $F$ is the homogeneous polynomial such that
$$
\split
 &F(x,y,1) =
     9 (x-1)y^9 
     - 3 (6 x^4 + 8 x^3 - 3 x^2 - 6 x + 1) y^6
 \\&\qquad
     + (9 x^7 + 24 x^6 + 13 x^5 - 8 x^4 - 11 x^3 - x^2 + 2 x - 1) y^3
     - (x+2)x^3(x^2-1)^3.
\endsplit
$$

(b). The degree $12$ curve $f(x,y,z)f(\omega x,y,z)f(\omega^2 x,y,z)=0$, where
$$
  f(x,y,z) = (x-y)(x+y)^3 + (2x+3y)z^3, \qquad \omega=e^{2\pi i/3},
$$
has $19$ ordinary triple points. The curve $f(x,y,z)=0$ admits parametrization
$t\mapsto(t^4-3t:t^4+2t:2t^3-1)$.

\endproclaim

In \S\sectSym\ we explain how the examples in Proposition~\propDegTen\ were constructed.

In \S\sectPerturb\ (which is independent of \S1) we obtain an irreducible curve of degree 10 with 12 triple points
as a perturbation of the dual Hesse arrangement (9 lines with 12 triple points)
with an additional line. Namely, we prove the following fact.

\proclaim{ Proposition \propPerturb }
There exists an analytic family  $C_u$, $|u|<1$, of curves of degree $10$ in $\CP^2$
such that $C_u$ for $u\ne 0$ is irreducible and has $12$ ordinary triple points, and
$C_0$ is given by the equation
$$
    (2x+2y-z)(x^3-y^3)(y^3-z^3)(z^3-x^3) = 0.                          \eqno(\eqTenLines)
$$
\endproclaim

We do not know if there exists an equisingular deformation of the curve
in Proposition \propDegTen(a) to those in Proposition \propPerturb.

\medskip\noindent
{\bf Remark \remGenLine.} Proposition \propPerturb\ holds for a generic line instead of the line $2x+2y=z$
in (\eqTenLines),
i.e., for any element of some Zariski open subset of $\check\P^2$ (see Remark~\remGLbis\ below).

\if01{
\medskip\noindent
{\bf Remark \remMethod.}
There are results on the existence of perturbations with prescribed
singularities in rather general settings (see e.g.~[\refGLSa], [\refGLSb], [\refS]).
Typically they deal with Cartesian equations of curves.
In our situation these results are not applicable, and we use another approach:
instead of perturbing a Cartesian equation we construct
a perturbation in parametric form. We do not know if non-equigeneric perturbations
(i.e. increasing the number of local branches) were studied somewhere
in a parametric form.
}\fi

\medskip\noindent
{\bf Remark \remMethod.}
There are results on the existence of perturbations with prescribed
singularities in rather general settings (see e.g.~[\refGLSa], [\refGLSb], [\refS]).
Those which deal with Cartesian equations are not applicable to our problem.
We construct the perturbation in a parametric form.
We do not know if non-equigeneric perturbations
(i.e. decreasing the number of irreducible components and/or increasing the sum of their genera)
were studied somewhere in a parametric form. It is interesting to see if our approach
could be used for such perturbations in more general settings.

\subhead Acknowledgement \endsubhead
I am grateful to Dmitry Kerner, Joaquim Ro\'e, and Eugenii Shustin for stimulating discussions.


\head\S\sectSym. Constructions based on a symmetry of order 3
\endhead
\subhead\sectBirat. Reduction by a birational transformation
\endsubhead
Let us show how the example from Proposition \propDegTen(a)² was found.
We look for a homogeneous polynomial $f(x,y,z)$ of degree $10$ of the form
$f(x,y,z) = f_1(x,y^3,z)$ such that the curve $C=\{f=0\}\subset\P^2$
is rational and has 12 triple points. Then the Newton polygon of 
$f_1(x,y,1)$ is contained in the trapezoid $[(0,0),(10,0),(1,3),(0,3)]$,
so, it is natural to consider the curve $C_1$ given by the equation $f_1(x,y,z)=0$
on the weighted projective plane $\P_{1,3,1}^2$ or on the Hirzebruch surface $\Sigma_3$,
which is $\P_{1,3,1}^2$ blown up at $(0:1:0)$. The curve $C$ is the preimage of $C_1$
under the three-fold branched covering $\rho:\P_{1,3,1}^2\to\P^2$, $(x:y:z)\mapsto(x:y^3:z)$.

It is easy to see that, under the above assumptions,
three triple points of $C$ are placed on the line $L=\{y=0\}$ so that their images
on $\Sigma_3$ are points where $C_1$ has a cubic tangency with the line $L_1=\rho(L)$.
The remaining nine triple points of $C$ are preimages of three triple points of $C_1$. 
The curves $L_1$ and $C_1$ (considered as curves on $\Sigma_3$) belong to the linear systems
$|G|$ and $|3G+F|$, where $F$ is a fiber and $G$ a general section of the projection $\Sigma_3\to\P^1$.
Let us blow up the three triple points of $C_1$ and then blow down the strict transforms of the fibers
passing through them. The resulting surface is $\Sigma_0=\P^1\times\P^1$. The strict transforms of $L_1$
and $C_1$ are smooth curves on it. We denote them by $L_2$ and $C_2$. They belong to the linear systems
$|A+3B|$ and $|3A+B|$ respectively, where $A$ 
is the strict transform of the $(-3)$-section
of $\Sigma_3$ and $B$ is the image of $F$.


\subhead\sectConstr. The construction
\endsubhead
Thus, in order to construct the desired curve of degree 10 in $\P^2$, it is enough to find smooth curves
$L_2\in|A+3B|$ and $C_2\in|3A+B|$ on $\Sigma_0$ such that $(L_2\cdot C_3)_{p_k}=3$, $k=1,2,3$, and to apply
the above transformations backward.
We observe that $L_2$ and $C_2$ are expressed
in $A$ and $B$ symmetrically. So, it is natural to search these curves in the form
$C_2=\{f_2(x,y)=0\}$ and $L_2=\{f_2(y,x)=0\}$ where $\deg_x f_2=1$ and $\deg_y f_2=3$
(here $(x,y)$ corresponds to $\big((x:1),(y:1)\big)\in\Sigma_0$).
We also observe that if $C_2$ is tangent to a line $x+y=2a$ at a diagonal point $p=(a,a)$, then
$(C_2\cdot L_2)_p=3$.

So, we set $f_2(x,y) = p(y)x+q(y)$ where $p=y^3-3y^2+1$ and $q=y^3-y^2+y$.
Then $C_2$ is tangent to the lines $x+y=2a$ for $a=0,\pm1$ at the
diagonal points $(0,0)$, $(1,1)$, $(-1,-1)$, thus it has a cubic tangency with $L_2$ at these points.
The passage from $\Sigma_0$ to $\Sigma_3$ should transform the equation $p(x)y+q(x)=0$ of $L_2$ into the equation $y=0$ of $L_1$.
This means that it transforms any equation $g(x,y)=0$, $\deg_y g=n$, into
$$
    g\Big(\frac{y-q(x)}{p(x)}\Big)p(x)^n = 0,
$$
and we obtain that the required curve $C$ in $\P^2$ is given by $f(x,y,z)=0$, where
$$
    f(x,y,1) = \left(x p\Big(\frac{y^3-q(x)}{p(x)}\Big) + q\Big(\frac{y^3-q(x)}{p(x)}\Big)\right)p(x)^3.
$$
In Figure \figSym\ the arrangement of $\R C_k$ and $\R L_k$ on $\R\Sigma_0$ and $\R\Sigma_3$ are shown
(we denote the real locus of a complex variety $X$ by $\R X$). The imaginary local branches of $C_1$
at triple points are schematically depicted by dashed lines; $\R\Sigma_k$ are represented by squares
whose opposite sides are glued according to the arrows.

\midinsert
\centerline{\epsfxsize=120mm\epsfbox{sym.eps}}
\botcaption{Figure \figSym}
    $\R L_2,\R C_2$ on $\R \Sigma_0$ and $\R L_1,\R C_1$ on $\R \Sigma_3$.
\endcaption
\endinsert

A parametrization can be computed as follows. It is clear that $C_2$ is parametrized by $y=s$, $x=-q(s)/p(s)$.
Then $C_1$ is parametrized by $s\mapsto(q(s):r(s):-p(s))\in\P_{1,3,1}^2$, where
$r(s) = 3(s-2)(s^2-1)^3s^3$. The parameter change $s=t^3+2$ yields a parametrization
$t\mapsto\big(q(s):r(s)^{1/3}:-p(s)\big)\big|_{s=t^3+2}$ of $C$.
The formulas in Proposition~\propDegTen\ are obtained by further
``cosmetic'' rescalings $(x,y,z)\to(x,3^{1/3}y,-z)$.


\subhead\sectQP. What about curves of degree 17 with 20 quadruple points
\endsubhead
This construction would be generalized to irreducible rational curves
of degree $n^2+1$ with $n^2+n$ ordinary $n$-fold points
as soon as there exist curves in the linear systems $|A+nB|$ and $|nA+B|$
on $\Sigma_0$ which have $n$ tangency points of order $n$. If we assume in
addition that the curves are symmetric to each other with respect to the diagonal
(as in \S\sectConstr), then the tangency points cannot be placed on the diagonal
when $n$ is even. We studied the symmetric case for $n=4$. In this case the
coordinates can be chosen so that the tangency points are at $(c,0)$, $(1,\infty)$
and at their symmetric images. Then, in the notation of \S\sectConstr,
$$
   f_2(x,y)=(y-1)(y^3+a_2 y^2 + a_1 y + a_0)x - (y-c)(y^3 + b_2 y^2 + b_1 y + a_0).
$$
The tangency conditions give a system of 6 polynomial equations
with 6 indeterminates. Using Gr\"obner bases, we checked that there are
no irreducible solutions.
In the general (i.e., asymmetric) case for $n=4$, the system of equations
is twice bigger and we did not succeed to solve it. We did not try $n\ge 5$.


\subhead\sectB. Three quartics with 19 triple points
\endsubhead
The construction of the curve in Proposition \propDegTen(b) is similar to the previous one.
We look for a curve of the form $F(x^3,y^3,z^3)=0$ which has a triple point at $p_1=(0:0:1)$,
three more triple points on each coordinate axis, and $9$ triple points somewhere else.
Then the quartic curve $F(x,y,z)=0$ passes through $p_1$, has a cubic tangency with each axis,
and has a triple point $p_2$ somewhere else. Let us apply the Cremona quadratic transformation
with the base points $p_1,p_2,p_3$, where $p_3$ is the point of cubic tangency with the axis $y=0$;
see Figure~\figCr\ on the left.
We obtain a cuspidal cubic $C_3$, a conic $C_2$, and four lines arranged as shown in Figure~\figCr\
on the right. If we forget the two lines passing through the point $q$ (see Figure~\figCr), the
remaining configuration is unique up to automorphism of $\P^2$, thus the whole picture is determined by
a choice of $q$ on $C_2$.
  If the coordinates are chosen as shown in Figure~\figCr,
then $C_3=\{x^3 = y^2 z\}$, $C_2 = \{ 2xy + xz = y^2 + 2yz\}$, and 
the equation in the statement of Proposition~\propDegTen(b) corresponds to $q=(3:2:8)$.

\midinsert
\centerline{\epsfxsize=120mm\epsfbox{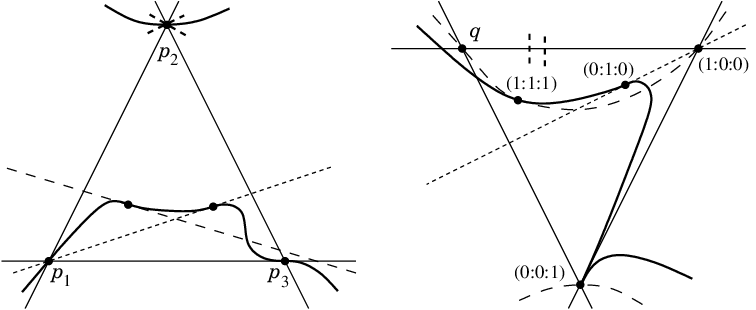}}
\botcaption{Figure \figCr} Cremona transformation in \S\sectB.
\endcaption
\endinsert

The parametrizations of the irreducible components of $F(x^3,y^3,z^3)=0$ can be found using the following observation:
if $t\mapsto\big(tx(t)^n:ty(t)^n:z(t)^n\big)$ is a parametrization of a curve $h(x,y,z)=0$, then the curve
$h(x^n,y^n,z^n)=0$ has $n$ irreducible components parametrized by
$$
    t\mapsto\big(\zeta^k tx(t^n):ty(t^n):z(t^n)\big),\qquad \zeta=e^{2\pi i/n}, \qquad k=1,\dots,n.
$$


\head\S\sectPerturb. Perturbation of the union of the dual Hesse arrangement with a generic line
\endhead
\subhead\sectDualHesse. The arrangement to perturb
\endsubhead
Denote the components of the dual Hesse line arrangement $(x^3-y^3)(y^3-z^3)(z^3-x^3)=0$ by
$$
\xalignat4
   &L_1=\{y=z\},            && L_4=\{z=x\},           && L_7=\{x=y\}, \\
   &L_2=\{y=\omega z\},     && L_5=\{z=\omega x\},    && L_8=\{x=\omega y\}, \\
   &L_3=\{y=\baromega z\},  && L_6=\{z=\baromega x\}, && L_9=\{x=\baromega y\},  && \omega = e^{2\pi i/3}.
\endxalignat
$$
It has twelve triple points $p_{ijk}=L_i\cap L_j\cap L_k$:
$$
\xalignat3
        &p_{123}=(1:0:0),              &     &p_{456}=(0:1:0),         &     &p_{789}=(0:0:1),\\
    p_1=&p_{147}=(1:1:1),              & p_4=&p_{159}=(\baromega:1:1), & p_7=&p_{249}=(1:\omega:1),\\
    p_2=&p_{258}=(1:\baromega:\omega), & p_5=&p_{267}=(1:1:\baromega), & p_8=&p_{168}=(\omega:1:1),\\
    p_3=&p_{369}=(1:\omega:\baromega), & p_6=&p_{348}=(1:\baromega:1), & p_9=&p_{357}=(1:1:\omega).
\endxalignat
$$
Let us denote $\Cal L_1=\{L_1,L_2,L_3\}$, $\Cal L_2=\{L_4,L_5,L_6\}$, $\Cal L_3=\{L_7,L_8,L_9\}$.

Let $L_0=\{2x+2y=z\}$. We parametrize it by $t\mapsto f_0(t)=(t:1-t:2)$.
Note that $L_0$ does not pass through the 12 triple points, nor through the intersections of the $L_i$'s
with the coordinate lines.
Let $t_1,\dots,t_9$ be such that $f_0(t_i)=L_0\cap L_i$, $i=1,\dots,9$, i.e.,
$$
   (t_1,\dots,t_9)=(-1,1-2\omega,1-2\baromega,\;\;\;2,2\baromega,2\omega,\;\;\;1/2,-\baromega,-\omega).    \eqno(\eqLzero)
$$
We are going to perturb $L_0\cup L_1\cup\dots\cup L_9$ so that the intersections of $L_0$ with the other lines are
smoothed while all the triple points are preserved.


\subhead\sectGenPerturb. A preliminary irreducible perturbation with double points
\endsubhead
Notice first that if $A$ is a rational curve in $\P^n$ and $L$ a line that
crosses $A$ at a point $p$ without tangency, then a smoothing of the intersection point can be described as follows
(cf.~[\refMO, Lemma~6.1]). Let us choose homogeneous coordinates $(x_0:\dots:x_n)$ in $\P^n$ so that
$L$ passes through $(1:0:\dots:0)$. Let $t\mapsto f(t)=(x_0(t):\dots:x_n(t))$ be a polynomial
parametrization of $A$ and let $t_0$ be such that $f(t_0)=p$. Let $C$ be the curve parametrized by
$t\mapsto(\hat x_0(t):\dots:\hat x_n(t))$ where
$\hat x_0(t)=(t-t_0-u)x_0(t)$ and $\hat x_i(t)=(t-t_0)x_i(t)$ for $i\ge 1$.
Then $C$ tends to $A\cup L$ when $u\to 0$. The degeneration of $C$ into $A\cup L$ and its analog
for pseudoholomorphic curves is often called {\it bubbling}.
This construction can be easily generalized to the case when several lines cross $A$ under condition that
each line passes through a base point of $\P^n$. Fortunately, this condition is satisfied by the dual Hesse arrangement.

So, we consider a curve $C=C_{u,\aa}$ whose parametrization
$$
     t\mapsto f(t)=f_{u,\aa}(t)=\big(tf_1(t) : (1-t)f_2(t) : 2f_3(t)\big),
$$
depends on $u\in\C$ and $\aa=(a_1,\dots,a_9;b_1,\dots,b_9)\in\C^9\times(\C\setminus\{0\})^9$
as follows:
$$
   f_\nu(t) = \prod_{i=1}^9 \big(t - t_i - (a_i + b_{i,\nu})u\big),
   \quad
   b_{i,\nu}=\cases b_i &\text{if $L_i\in\Cal L_\nu$,}\\
   0 &\text{otherwise.}\endcases
$$
If $\aa$ is fixed and $u$ tends to zero, then $C$ tends to $L_0\cup L_1\cup\dots\cup L_9$.
By construction, $C$ has triple points at $(1:0:0)$, $(0:1:0)$, and $(0:0:1)$.
However, each of the other nine triple points in general splits into three nodes.
Our aim is to find an analytic germ $u\mapsto\aa(u)$ such that 
all the 12 triple points are preserved in $C_{u,\aa(u)}$.


\subhead\sectStepOne. A perturbation with the triple points preserved up to $O(u^2)$
\endsubhead
From now on we pass to the affine chart $z\ne 0$, which we identify with $\C^2$.
So, we denote a point $(x:y:z)$ by $(x/z,y/z)$ and say that $C$ is parametrized by
$$
    t\mapsto f(t)=\big(X(t),Y(t)\big), \qquad
    X(t)=\frac{tf_1(t)}{2f_3(t)}, \qquad Y(t)=\frac{(1-t)f_2(t)}{2f_3(t)},
$$
where $f_\nu(t)$ are as above. Then
$$
  X = \frac{t}{2}\frac{A_1 A_2 A_3\, B_7 B_8 B_9}
                         {B_1 B_2 B_3\, A_7 A_8 A_9},
  \quad
  Y = \frac{(1-t)}{2}\frac{A_1 A_2 A_3\, B_4 B_5 B_6}
                           {B_1 B_2 B_3\, A_4 A_5 A_6},
  \quad
  \cases A_i = t - t_i - a_i u,\\ B_i = A_i - b_i u.\endcases
$$

If $p_k\in L_i$ and $L_i\in\Cal L_j$, we set $i_j(k) = i$.
Let $\tau_{k,j}=\tau_{k,j}(\aa)$ be found from the condition
$$
     \lim_{u\to 0} f(t_i + \tau_{k,j}u) = p_k, \qquad i=i_j(k).
$$
Each $\tau_{k,j}$ is of the form $a_i+\lambda_{k,j} b_i$, $i=i_j(k)$, $\lambda_{k,i}\in\Q(\omega)$.
For example, we have $p_1=(1,1)$, $i_2(1)=4$, $t_4=1/2$  (see (\eqLzero)), and
$$
  \lim_{u\to 0} Y(t_4+\tau u) = \frac{(1-t_4)(\tau-a_4-b_4)}{2(\tau-a_4)},
$$
hence $\tau_{1,2}=a_4+\tfrac{1}{3}b_4$. Let us set
$$
      p_{k,j}(u) = \cases f(t_i + \tau_{k,j}u), & u\ne 0,\; i=i_j(k),\\ p_k, & u=0. \endcases
$$
Both components of the derivatives $p'_{k,j}(0)$ are linear combinations of $a_1,\dots,a_9$, $b_1,\dots,b_9$ with
coefficients in $\Q(\omega)$, for example, the $X$-component of $6 p'_{1,2}(0)$ is
$$
   -2 b_1 + (2+4\omega)(b_2-b_3) + 3 a_4 + b_4 + 4 b_7 + (4+2\omega)b_8 + (2-2\omega)b_9.
$$

Let $k\in\{1,\dots,9\}$.
We are going to find a condition on $\aa$ to assure that the triple point at $p_k$ is preserved in $C_{u,\aa}$ up to $O(u^2)$
(i.e., the distances between the double points of $C_{u,\aa}$ appearing near $p_k$
are of that order). To this end, in a neighborhood of $p_k$,
we approximate $C_{u,\aa}$ by
$L'_1\cup L'_2\cup L'_3$, where $L'_j$ is the line parallel to $L_{i_j(k)}$ passing through
the point $p_k + u p'_{k,j}(0)$. These lines are concurrent if and only if there exist $\xi_1,\xi_2,\xi_3$ such that
$$
     v_1 + \xi_1 w_1 = v_2 + \xi_2 w_2 = v_3 + \xi_3 w_3,
$$
where
$v_j = p'_{k,j}(0)$ and $w_1=(1,0)$, $w_2=(0,1)$, $w_3=p_k$ (the vector $w_j$ is parallel to $L_{i_j(k)}$).
This is a system of four linear non-homogeneous equations for three unknowns,
thus it has a solution if and only if $e_k=0$, where
$$
    e_k = \det\left(\matrix w_1  &  -w_2  &    0    &  v_2 - v_1 \\
                            w_1  &    0   &  -w_3   &  v_3 - v_1 \endmatrix\right)      \eqno(\eqDet)
$$
(here the $v_j$ and $w_j$ are interpreted as columns). Since the $w_j$ are constant,
$e_k$ is a linear combination of $a_1,\dots,a_9,b_1,\dots,b_9$ with coefficients in $\Q(\omega)$,
for example, we have
$$
   e_1 = \tfrac12(a_1 + b_1 + a_4 + b_4) - 4(a_7+b_7) + b_8+b_9 + \tfrac{2\omega-1}{7}(b_2+b_6) + \tfrac{2\baromega-1}{7}(b_3+b_5).
$$
Thus we obtain a system of homogeneous linear equations $e_1=\dots=e_9=0$.
Each solution with all $b_i$'s non-zero yields a desired perturbation with all the triple points
preserved up to $O(u^2)$. For such a solution we can choose
$$
  \wh\aa=\left(-\tfrac{1}{3},\alpha,\bar\alpha,
              -\tfrac{1}{3},\bar\alpha,\alpha,
               \tfrac{37}{42},0,0;\;\;
               1,1,1,1,1,1,-\tfrac{1}{2},1,1\right),
               \quad
               \alpha=\frac{38\omega-5}{21}.                          \eqno(\eqA)
$$

\noindent
{\bf Remark \remUseSym.} If $L_0$ is generic, then the existence of a solution
with all $b_i$'s non-zero can be derived from the symmetry of the Hesse configuration. Indeed,
suppose that for any solution there exists $i$ such that $b_i=0$ in this solution.
Then there exists $i$ such that $b_i=0$ in any solution.
Since the configuration is symmetric
and $L_0$ is generic, it follows that $b_1=\dots=b_9=0$ in any solution.
This fact would imply that the coefficient of each $a_i$ in each $e_k$ is zero,
which is not the case.

However, we do not know how to do the next step (in \S\sectStepTwo)
without the above computations for a concrete line
(for which we have chosen $2x+2y=z$ though this particular choice is not at all important).


\subhead\sectStepTwo. A perturbation with all triple points preserved
\endsubhead
Roughly speaking, we just blow up the space of parameters at $(0,\wh\aa)$ (see (\eqA)) and apply the implicit
function theorem. Let us explain it in more detail.

\midinsert
\centerline{\epsfxsize=75mm\epsfbox{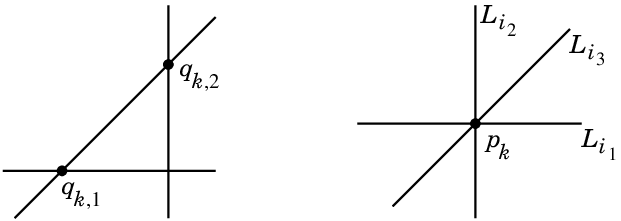}}
\botcaption{Figure \figPerturb}
    The double points of $C_{u,\aa}$ that converge to $p_k$.
\endcaption
\endinsert

For a fixed $k\in\{1,\dots,9\}$, let $q_{k,j}(u,\aa)$, $j=1,2$, be the double point of $C_{u,\aa}$
whose local branches converge to the germs of $L_{i_j(k)}$ and $L_{i_3(k)}$ at $p_k$ as $u\to 0$ (see Figure \figPerturb).
Let $s_{k,1}(u,\aa)$ and  $s_{k,2}(u,\aa)$ be such that
$$
     f_{u,\aa}\big(t_{i_3(k)} + s_{k,j}(u,\aa)u\big) = q_{k,j}(u,\aa), \qquad j=1,2.
$$
i.e. $t_{i_3(k)} + s_{k,j}(u,\aa)u$ is the parameter of the center of the local branch of $C_{u,\aa}$ at $q_{k,j}(u,\aa)$
that converges to $(L_{i_3(k)},p_k)$.
The condition that the triple point $p_k$ is preserved under the
perturbation is equivalent to the condition that
$$
     s_{k,1}(u,\aa)=s_{k,2}(u,\aa).                        \eqno(\eqPhi)
$$

Let $\wh\aa$ be as in (\eqA). The fact that it is a solution of the equations $e_1=\dots=e_9=0$ (see \S\sectStepOne)
implies that (\eqPhi) holds up to $O(u^2)$ when $\aa = \wh\aa+\aaa u$, $\aaa\in\C^{18}$, i.e.,
$$
    s_{k,2}(u,\wh\aa+\aaa u) - s_{k,1}(u,\wh\aa+\aaa u) = u^2 \Phi_k(u,\aaa),
$$
where $\Phi_k$ is analytic in a neighborhood of $\{(0,\aaa)\mid\aaa\in\C^{18}\}$ in $\C\times\C^{18}$.

One can compute that $\Phi_k(0,\aaa)$ are affine linear functions of $\aaa$ with coefficients in $\Q(\omega)$
(see more details in the recursion step for $n=1$ in \S\sectSeries),
for example, we have $\Phi_k(0,\aaa)=c_ke_k(\aaa)+d_k$ with $c_k,d_k\in\Q(\omega)$ for $k=1,5,9$
though the linear part of $\Phi_k(0,\aaa)$ is not proportional to $e_k(\aaa)$
for other values of $k$ (recall that $e_k$ is defined in (\eqDet)).
Moreover $\Phi|_{u=0}$ is surjective for
$$
    \Phi=(\Phi_1,\dots,\Phi_9):\C\times\C^{18}\to\C^9,
$$
in particular, $\det\big(\Phi|_{0\times E}\big) = 64 / (3^{12}\times 49)$, where
$$
    E = \{\aaa\in\C^{18}\mid \tilde a_1=\tilde a_4,\;\;
    \tilde a_7=\tilde b_1 = \tilde b_4 = \tilde b_5 = \tilde b_6 = \tilde b_7 = \tilde b_8 = \tilde b_9 = 0\}. \eqno(\eqE)
$$

Hence, by the implicit function theorem, for a neighborhood $U$ of $0$ in $\C$
there exists an analytic mapping $U\to E$, $u\mapsto\aaa(u)$, such that $\Phi(u,\aaa(u))=0$.
Then $t\mapsto f_{u,\wh\aa+u\aaa(u)}(t)$ is the required perturbation of $L_0\cup L_1\cup\dots\cup L_9$
into a rational curve with 12 triple points.

\medskip\noindent
{\bf Remark \remGLbis.} (Cf.~Remarks \remGenLine\ and \remUseSym.) Proposition~\propPerturb\ remains true for
a generic line $L_0$ because the set of lines such that $\text{rank}(\Phi|_{u=0})=9$ is Zariski open and non-empty.

\medskip\noindent
{\bf Remark \remHDbis.}
One might think that the same arguments would yield an irreducible
perturbation of a curve $(ax+by+cz)(x^n-y^n)(y^n-z^n)(z^n-x^n)$ with $n^2$ triple points and three $n$-fold points
for $n=4$ and $5$. The results of \S\sectStepOne\ indeed extend to these cases. However, the analog
of $\Phi|_{u=0}$ is no longer surjective: it is a mapping $\C^{24}\to\C^{16}$ of rank $15$ for $n=4$, and
$\C^{30}\to\C^{25}$ of rank $20$ for $n=5$.


\subhead\sectSeries. Power series expansions
\endsubhead
Let $\aaa(u)$ be as in \S\sectStepTwo. Denote $\aa(u)=\wh\aa+u\aaa(u)$, $C_u=C_{u,\aa(u)}$
and $f_u=f_{u,\aa(u)}$. Let $p_k(u)$ be the triple point of $C_u$ which
tends to $p_k$ as $u\to 0$, and let $T_{k,j}(u)\in f_u^{-1}(p_k(u))$ be such that
the germ of $f_u$ at $T_{k,j}(u)$ defines the local branch of $C_u$ at $p_k(u)$ that tends to 
$(L_{i_j(k)},p_k)$ as $u\to 0$. By the definition of $\tau_{k,j}$ (see \S\sectStepOne)
there exist analytic functions $S_{k,j}(u)$ such that
$$
    T_{k,j}(u) = t_{i_j(k)} + u\tau_{k,j}\big(\aa(u)\big) + u^2 S_{k,j}(u).
$$

In this subsection we explain how to recursively compute the power series expansion of $\aa(u)$
simultaneously with the expansions of the $S_{k,j}(u)$. Let us set:
$$
\split
    &\aa(u) = \sum_{n=0}^\infty \aa_n u^n,
    \quad \aa^{[n]}(u)=\sum_{m=0}^n\aa_m u^m,
    \quad \aa_n=(a_{1,n},\dots,a_{9,n};\,b_{1,n},\dots,b_{9,n}),
    \\
    &S_{k,j}(u) = \sum_{n=0}^\infty S_{k,j,n}u^n,
    \qquad S_{k,j}^{[n]}(u) = \sum_{m=0}^n S_{k,j,m} u^m.
\endsplit
$$
Then $C_{u,\aa^{[n]}(u)}$ is a perturbation with all triple points preserved up to $O(u^{n+2})$.

The recursive computation of $\aa^{[n]}$ and $S_{k,j}^{[n]}$ goes as follows. 
Note that all the equations below are linear (because we expand $S_{k,j}$ instead of $T_{k,j}$).
Moreover, the homogeneous part of (\eqSys) does not depend on $n$.

\medskip\noindent
$1^\circ$. {\sl Initialization of the recursion.}
For each $k=1,\dots,9$, find $s_1,s_2,s_3$ from the system of equations
$$
      f_{u,\wh\aa}\big(t_{i_j(k)}+u\tau_{k,j}(\wh\aa)+u^2 s_j\big)=
      f_{u,\wh\aa}\big(t_{i_3(k)}+u\tau_{k,j}(\wh\aa)+u^2 s_3\big)+O(u^2), \quad j=1,2,
$$
(these are four linear equations for three unknowns,
but a solution exists due to the definition of $\wh\aa$).
Then set $S_{k,j}^{[0]}=s_j$, $j=1,2,3$, and $\aa^{[0]}=\wh\aa$.
\medskip\noindent
$2^\circ$. {\sl The recursion step.}
Assuming $\aa^{[n-1]}(u)$ and all the $S_{k,j}^{[n-1]}(u)$ known, compute 
 $\aa^{[n]}(u)$ and all the $S_{k,j}^{[n]}(u)$ as follows.
\roster
\item"$\bullet$"
For each $k=1,\dots,9$ and $j=1,2$, find $s_{k,j}=s_{k,j}(\aa_n)$ and $s^*_{k,j}=s^*_{k,j}(\aa_n)$ from the equation
$$
   \phi_{k,3}^{[n]}(0,s_{k,j}) = \phi_{k,j}^{[n]}(0,s^*_{k,j}),
$$
where $\phi_{k,j}^{[n]}$ is the analytic function such that
$$
   f_{u,\aa^{[n]}(u)}\Big(t_{i_j(k)}+u\tau_{k,j}\big(\aa^{[n]}(u)\big)
                       +u^2\big(S_{k,j}^{[n-1]}(u)+s u^n\big)\Big) = u^{n+1}\phi_{k,j}^{[n]}(u,s).
$$
\item"$\bullet$"
Find $\aa_n\in E$ (see (\eqE)) from the system of equations
$$
           s_{k,1}(\aa_n)=s_{k,2}(\aa_n), \qquad k=1,\dots,9.                      \eqno(\eqSys)
$$
\item"$\bullet$"
Set
$S_{k,j,n} = s_{k,j}^*(\aa_n)$, $j=1,2$, and
$S_{k,3,n} = s_{k,1}(\aa_n)$.
\endroster
The result of computation of $\aa^{[7]}$ and $S_{j,k}^{[7]}$, as well as the Maple code
with comments, is available at
{\tt\ https://www.math.univ-toulouse.fr/\~{}orevkov/12tp.html}.


\enddocument